\documentclass{article}
\usepackage[utf8]{inputenc}
\usepackage{amssymb}
\usepackage{amsthm}
\usepackage[english]{babel}
\usepackage{mathtools}
\usepackage{xcolor}
\usepackage{comment}
\usepackage{bbm}
\usepackage{cancel}
\usepackage{setspace}
\usepackage{multicol}
\usepackage{textcomp}

\newtheorem{thm}[equation]{Theorem}
\numberwithin{equation}{section}
\newtheorem{lem}[equation]{Lemma}
\newtheorem{prop}[equation]{Proposition}
\newtheorem{rem}[equation]{Remark}
\newtheorem{cor}[equation]{Corollary}
\newtheorem{defn}[equation]{Definition}

\newcommand\Z
 {\mathbb{Z}}

 \newcommand\F
 {\mathbb{F}}
 
\newcommand\fs
 {\mathcal{F}}

\newcommand\1
 {\mathbbm{1}}

\newcommand\into 
{\hookrightarrow}

\usepackage{tikz-cd}

\usepackage[margin=2.5cm]{geometry}
\title{On the modular cohomology of $GL(2,p^n)$ and $SL(2,p^n)$}
\author{}
\date{ }

\begin{document}
 \begin{center}
     \huge{On the modular cohomology of $GL(2,p^n)$ and $SL(2,p^n)$}\\
     
\vspace{0.3cm}

     \normalsize{Anja Meyer}
 \end{center}

\vspace{0.4cm}

\begin{center}
\textbf{Abstract}\\
    \small{Let $p$ be an odd prime. Denote a Sylow $p$-subgroup of $GL_2(\Z/p^n)$ and $SL_2(\Z/p^n)$ by $S_p(n,GL)$ and $S_p(n,SL)$ respectively. The theory of stable elements tells us that the mod-$p$ cohomology of a finite group is given by the stable elements of the mod-$p$ cohomology of it's Sylow $p$-subgroup. We prove that for suitable group extensions of $S_p(n,GL)$ and $S_p(n,SL)$ the $E_2$-page of the Lyndon-Hochschild-Serre spectral sequence associated to these extensions does not depend on $n>1$. Finally, we use the theory of fusion systems to describe the ring of stable elements.}
\end{center}

\vspace{0.4cm}

\tableofcontents

\section{Introduction}

Cohomology of matrix groups is a well-studied and widely applicable area of mathematics, whose roots date back to the very beginnings of (co-)homology theories about a hundred years ago. Much work has been done over the century. In Quillen's famous 1972 paper \textit{On the cohomology and K-theory of the general linear group over a field} [Q] he established $H^*(GL_k(\Z/p, \F_q)$ for $p,q$ distinct primes and any $k$. The cohomology of $GL_2(\Z/p)$ with mod-$p$ coefficients was calculated by Aguade in 1980 in [A]. In 2020 Diaz, Garaialde, Mazza and Park moved to the profinite case and computed the mod-$p$ cohomology of $GL_2(\Z_p)$ in [DGMP] (published in 2023). Despite a vast amount of work by various mathematicians, the mod-$p$ cohomology of $GL_m(\Z/p^n)$ and $SL_m(\Z/p^n)$ for $n>1$ and $m>1$ remains unknown. In this paper we introduce a method to tackle this open problem for $m=2$ and odd primes $p$ and present results that can be used for these calculations. The aim in the long run is to generalise the here presented method to arbitary $m$.\\
The paper is organised as follows: In Chapter 2 we examine $GL_2(\Z/p^n)$ and $SL_2(\Z/p^n)$ and some subgroups with group theoretic tools. Let $S_p(n,GL) \in Syl_p(GL_2(\Z/p^n))$ and $S_p(n,SL) \in Syl_p(SL_2(\Z/p^n))$. In Chapter 2 we will set up short exact seqeuences
\begin{center}
    $1 \rightarrow L_n \rightarrow S_p(n,GL) \rightarrow S_p(1,GL) \rightarrow 1$ 
      \hspace{0.5cm} and  \hspace{0.5cm}
    $1 \rightarrow K_n \rightarrow S_p(n,SL) \rightarrow S_p(1,SL) \rightarrow 1$\\
\end{center}
where $L_n$ and $K_n$ are the respective kernels of reduction modulo $p$.\\
In Chapter 3 we examine $L_n$ and $K_n$ through the lense of profinite group theory and gather information their modular cohomology. Theorem \ref{Cohom(Kn)} gives an explicit formulation.\\
This will be used in Chapter 4 where we set up the Lyndon-Hochschild-Serre spectral sequence. Theorem \ref{G-module isomorphism} is one of the two main theorems of this paper. It tells us that the $E_2$ pages associated to the group extensions above do not depend on $n$. This theorem makes this work here possible.\\
In Chapter 5 we look into the stable elements. The theory of fusion systems is needed, and we will recall the required notions before stating and proving the main theorem of this paper. This is Theorem \ref{final theorem} and gives and explicit description of the modular cohomology rings of $GL_2(\Z/p^n)$ and $SL_2(\Z/p^n)$.\\

\noindent \textbf{Acknowledgements:} Chapters 2-4 are part of the author's PhD work, and she would like to express her gratitude to Peter Symonds for his supervision and guidance. In chapter 5 the idea to restrict our attention to the centric radical subgroups was proposed by Ran Levi and the author is grateful for his thoughts. Furthermore, she wishes to thank Adam Jones for helpful comments on a previous draft of this paper.

\section{Definitions and basic properties}\label{section Kn}

\noindent \textbf{Notation:}  $\1$ denotes the $2 \times 2$-identity matrix\\

\noindent In order to understand the groups in questions and some of their suitable subgroups, we start by considering reduction modulo $p^m$ for $m<n$. The kernel of this map will turn out to be a key player in this paper.\\

\noindent For $n>m>0 \in \Z$, define $\gamma_{n,m}:GL_2(\Z/p^n ) \rightarrow GL_2(\Z/p^m )$ to be reduction modulo $p^m$. 
The kernel of $\gamma_{n,m}:GL_2(\Z/p^n ) \rightarrow GL_2(\Z/p^m )$ is
\begin{center}
        $L_{n,m}=\{\1+p^mX|X\in M_{2\times 2}(\Z/p^n)\} \triangleleft GL_2(\Z/p^n),$\\
 
    \end{center}
Furthermore, the restriction of $\gamma_{n,m}$ to $ SL_2(\Z/p^n) \rightarrow SL_2(\Z/p^m)$ has kernel 
  \begin{center}
        $K_{n,m}=\{\1+p^mX|X\in M_{2\times 2}(\Z/p^n\Z), \ \textnormal{det}(\1+p^mX)=1\} \triangleleft SL_2(\Z/p^n ).$
    \end{center}

\noindent One verifies that $\gamma_{n,m}$ is surjective and arrives at the following short exact sequences
\begin{center}
    $ 1 \rightarrow L_{n,m} \into GL_2(\Z/p^n) \xrightarrow{\gamma_{n,m}} GL_2(\Z/p^m) \rightarrow 1$,

    $ 1 \rightarrow  K_{n,m} \into SL_2(\Z/p^n) \xrightarrow{\gamma_{n,m}} SL_2(\Z/p^m) \rightarrow 1 $.
\end{center}

We examine some properties of $L_{n,m}$ and $K_{n,m}$. We start with the case $m=n-1$. 

\begin{prop}\label{elem.ab.} For $n>1$ the groups $L_{n,n-1}$ and $K_{n,n-1}$ are elementary abelian. In particular,\\
    $ L_{n,n-1} \cong (\Z/p )^4$ and $K_{n,n-1} \cong (\Z/p )^3$. Furthermore, for $\1+p^{m-1}X \in L_{m,m-1}$ we can choose the explicit isomorphism
    \begin{center}
        $\theta:L_{m,m-1} \rightarrow L_{n,n-1}$\\
        $X \mapsto p^{n-m}X$.
    \end{center}
    Note that this map restricts to $K_{m,m-1}$, too. 
\end{prop}
 \textbf{Proof:} Since $p^k$ generates a subgroup in $\Z/p^{k+1}$ of order $p$, $|L_{n,n-1}|=p^4$ as there are four positions $a,b,c,d$ to fill in the matrix  
     $A=\1+p^{n-1}X=\1+p^{n-1}\begin{pmatrix} a & b \\ c & d  \end{pmatrix}$.\\
 Furthermore for $A,B \in L_{n,n-1}$,
\begin{align*}
     AB &=\1+p^{n-1}X+p^{n-1}Y+p^{2(n-1)}XY \\
        &\equiv \1+p^{n-1}(X+Y) \mod p^n. 
\end{align*}

Since matrices commute under addition, we see that $L_{n,n-1}$ is abelian. Next we show that $\1+p^{n-1}X\in L_{n,n-1}$ has order $p$. \\
As $A^2 \equiv \1+p^{n-1}(X+X) \mod p^n$,
\begin{align*}
     A^p  &\equiv \1+p^{n-1}(pA) \mod p^n\\
      &\equiv \1 \mod p^n.
\end{align*}

By the Structure Theorem $L_{n,n-1}$ is isomorphic to $(\Z/p)^4$.\\
Next, we turn our attention to
  \begin{center}
        $\theta:L_{m,m-1} \rightarrow L_{n,n-1}$\\
        $X \mapsto p^{n-m}X$.
    \end{center}
Take $A=\1+p^{n-1}Y \in L_{n,n-1}$, $Y\in M_{2\times 2}(\Z/p^n).$ As we reduce modulo $p^n$ in $L_{n,n-1}$ we can reduce $Y$ to $Y'\in M_{2 \times 2}(\Z/p)$ and $\1+p^{n-1}Y =\1+p^{n-1}Y' $. Let now $B=\1+p^{m-1}Y'$, then $B\in L_{m,m-1}$. Applying the map $\theta$ gives \begin{center}
    $\theta(B)=\theta(\1+p^{m-1}Y')=\1+p^{m-1}(p^{n-m}Y')=1+p^{n-1}Y'=A$
\end{center} as required. As inverse to $\theta$ we simply send
$1+p^{n-1}X$ to $1+p^{m-1}X$, which can be checked to be a group homomorphism, too.\\
For $K_n$ we need to adjust the above slightly. As $K_{n,n-1}\subset SL_2(\Z/p^n)$, for any $A\in K_{n,n-1}, \ \textnormal{det}(A)=1$. A computation shows that this implies that $tr(A)=2$, which imposes a linear relation on $K_{n,n-1}$ resulting in $|K_{n,n-1}|=p^3$ rather than being equal to $p^4$. Since $K_{n,n-1}$ is abelian as well, we use the structure theorem and see that $K_{n,n-1}\cong (\Z/p)^3$.\\
Finally we turn our attention to $\theta|_{K_{m,m-1}}$. Following the argument further up in the proof shows $\theta|_{K_{m,m-1}}$ to be an isomorphism, too. \hfill $\square$\\

\noindent \textbf{Notation:} Whenever $m=1$ we write $L_n$ and $K_n$ to simply notation.

\begin{prop} The group $L_n$ is not abelian for $n>2$. Likewise, $K_n$  not abelian for $n>2$.
\end{prop}
\textbf{Proof:} We consider the following matrices in $K_n$: $A=\begin{pmatrix} 1  & 0 \\ p & 1 \end{pmatrix}, \ B= \begin{pmatrix} 1 & p \\ 0 & 1 \end{pmatrix}$. Multiplying these shows that $AB$ and $BA$ are not congruent mod $p^n$ for $n>2$. Hence $K_n$ is not abelian for $n>2$. Since $K_n$ is a subgroup of $L_n$ it follows that $L_n$ is not abelian for $n>2$ either. \hfill $\square$\\

\begin{lem}\label{p groups}
The groups $L_n$ and $K_n$ are p-groups.
\end{lem}
\textbf{Proof:}
Recall that by Proposition \ref{elem.ab.} $|L_{n,n-1}|=p^4$ and consider the short exact sequence
\begin{center}
    $1 \rightarrow L_{n,n-1}\rightarrow L_n \rightarrow L_{n-1}\rightarrow 1 $
\end{center}
from which follows that $|L_n|=|L_{n-1}| \times p^4$. By induction on $n$ we get that $|L_n|=p^{(n-1)4}$, making $L_n$ a p-group. \\
For $K_n$ we recall that $|K_n|=p^3$ and applying the above argument yields $|K_n|=p^{(n-1)3}$. Hence $K_n$ is a p-group as well.

\hfill $\square$\\

Recall that we wish to find suitable group extensions for the Sylow $p$-subgroups of $GL_2(\Z/p^n )$ and $SL_2(\Z/p^n)$, in order to compute their cohomology. We will now set up these extensions, serving as preparation for Section 4.

\begin{thm}\label{def Sylow} The set of matrices of the form\[\begin{pmatrix}
1+pa & b\\
pc & 1+pd
\end{pmatrix} \in GL_2(\Z/p^n )\] 
where $a,b,c,d\in \Z/p^n $ form a Sylow $p$-subgroup of $ GL_2(\Z/p^n)$. Denote this group $S_p(n,GL)$.

The set of matrices of the form\[\begin{pmatrix}
1+pa & b\\
pc & 1+pd
\end{pmatrix} \in SL_2(\Z/p^n)\] 
where $a,b,c,d\in \Z/p^n $ with determinant equal to 1 form a Sylow $p$-subgroup of $ SL_2(\Z/p^n )$. Denote this group $S_p(n,SL)$.

\end{thm}

\textbf{Proof:} First note that $GL_2(\Z/p)$ and $SL_2(\Z/p )$ Sylow-$p$-subgroups are isomorphic to $C_p$ - the cyclic group on $p$ elements. We recall the following elmentary fact:\\
Let $G$ be a finite group and $N \triangleleft G$ be a $p$-group with order $p^n$ for some $n$, and let $S\in Syl_p(G)$. Under a projection $\pi:G \rightarrow G/N$, $\pi (S)$ is a Sylow-$p$-subgroup of $G/N$. Focussing on $GL_2(\Z/p^n)$ we apply this fact to
\begin{center}
$1 \rightarrow L_n \rightarrow GL_2(\mathbb{Z}/p^n ) \xrightarrow{\gamma_{n,1}} GL_2(\mathbb{Z}/p) \rightarrow 1, $
\end{center}
where we note that $GL_2(\Z/p^n)/L_n =GL_2(\mathbb{Z}/p )$. A Sylow $p$-subgroup of $GL_2(\mathbb{Z}/p )$ is a Sylow $p$-subgroup of $\gamma(GL_2(\mathbb{Z}/p^n))$. Let $a,b,c,d \in \Z/p^n$ the inverse image of an element in $S_p(1,G)$ is 
\[
\begin{pmatrix}
1+pa & b\\
pc & 1+pd
\end{pmatrix}.
\] 
The group of matrices of this form has order $p^{4n-3}$. As we have 
\[
|GL_2(\Z/p^n)|=|GL_2(\Z/p )|\times p^{4(n-1)}=p^{4n-3}(p-1)^2(p+1),
\] 
we see that the order of the group above divides the order of $GL_2(\mathbb{Z}/p^n )$ and is thus a Sylow $p$-subgroup of $GL_2(\mathbb{Z}/p^n)$. \\

\noindent The proof for $S_p(n,S)$ is analogous but with the requirement that \[
\textnormal{det}\begin{pmatrix}
1+pa & b\\
pc & 1+pd
\end{pmatrix}=1.
\] 
Thus the group of such matrices has order $p^{3n-2}$. As \[
|SL_2(\mathbb{Z}/p^n )|=|SL_2(\Z/p)|\times p^{3(n-1)}=p^{3n-2}(p-1)(p+1).
\] 
As the highest power of $p$ is $p^{3n-2}$, this group is indeed a Sylow $p$-subgroup of $SL_2(\Z/p^n\Z)$. \hfill $\square$\\

\vspace{0.3cm}
 \begin{rem}
     The above theorem gives that we have short exact sequences
     \begin{center}
    $ 1 \rightarrow L_n \rightarrow S_p(n,GL) \rightarrow S_p(1,GL)\cong C_p \rightarrow 1 \hspace{2cm}$ \\
    $ 1 \rightarrow K_n \rightarrow S_p(n,SL) \rightarrow S_p(1,SL) \cong C_p \rightarrow 1\hspace{2cm}$. 
\end{center}

 \end{rem}

\noindent These will be used in the Lyndon-Hochschild-Serre spectral sequence in Chapter 4. In order to use these we will need to know the modular cohomology of $L_n$ and $K_n$ respectively for any $n$. As the order of $S_p(n,GL)$ and $S_p(n,SL)$ grows rapidly as $n$ gets larger, so does the number of their subgroups. Classical group theory does not give us sufficient information. Instead we will turn to the theory of pro-finite group. This study is presented in the next chapter.

\section{Profinite methods}

We recall the following definitions.

\begin{defn}[Profinite Group]
Let $(G_i, \pi_m^n)$ be an inverse system of finite topological groups and form the group $ \prod_i G_i$ for $i \in \Z$.  We collect at all the elements 
\begin{center}
    $\left\{ g_i \in \prod_i G_i| \text{for  } m\leq n, \ \pi_m^n (g_n) = g_m \right\}$.
\end{center}  
These elements form a subgroup $G$ of $\prod_i G_i$, with the induced topology.  This subgroup is called a \textit{profinite group} and is denoted by 
\begin{center}
    $G=\underleftarrow{\lim} G_i$.
\end{center} 
The map $\pi_n:G \rightarrow G_n$ is the canonical projection from G onto the $n^{th}$ group in the inverse system. 
\end{defn}
Any finite group is an example of a profinite group with a constant inverse system. A profinite group $G$ is called a pro-$p$-group if all $G_i$ are $p$-groups. Our next aim is building the theory to quote and then apply a result from [MS], which will be Theorem \ref{mis} in this paper, to determine the mod-$p$ cohomology of $K_n$ and $L_n$. We require the following technical background from profinite group theory:

\begin{defn}[Powerfully embedded]
A closed normal and finitely generated subgroup N of a pro-p-group G is powerfully embedded in G if $[G,N]\subset N^{2p}$
\end{defn}

\begin{defn}[Powerful group]
A pro-p-group group G is powerful if it is
powerfully embedded in itself.
\end{defn}

We continue with the following technical result:

\begin{prop}\label{powerful}
For $p$ an odd prime $K_n$ and $L_n$ are powerful pro-$p$-groups for all $n>1$.
\end{prop}
\noindent \textbf{Proof:} We start by proving the theorem for $K_n$ first. As $K_n \subset SL_2(\Z/p^n)$ we have the added restriction of the determinat of the matrices being equal to 1. The proof for $L_n$ will then follow.\\
Note that by Lemma \ref{p groups}, $K_n$ is a $p$-group. This, with respect to the constant inverse system, makes it a pro-$p$-group. 

\noindent To show that $K_n$ is powerful we note the following implications from [MS]: 
\begin{center}
    $\big( [G,G]\subseteq G^p\big) \iff \ \big(G$ almost powerfully embedded in itself$\big) \ \implies \ \big(G$ powerful$\big)$,
    
\end{center}
where $G^p$ is the closure of the subgroup of G generated by $p^{th}$ powers of elements of $G$. If $G$ is a finite group, $G^p$ is simply the subgroup of $G$ generated by $p^{th}$ powers of elements in $G$.\\

\noindent 
We show that $[K_n,K_n]\subseteq K_n^p$ by induction on n. For the base case we take $K_2=K_{2,1}$ which we have already proven to be elementary abelian. Hence $K_2^p \cong \{e\}$, the trivial subgroup and $[K_2,K_2]\subseteq K_2^p$. We assume $[K_{n-1},K_{n-1}]\subseteq K_{n-1}^p$. Then $[K_n,K_n]\subseteq K_n^p\cdot K_{n,n-1}$ which we see by using the short exact sequence 
\begin{center}
    $1 \rightarrow K_{n,n-1}\rightarrow K_n \rightarrow K_{n-1} \rightarrow 1$,
\end{center} giving $K_{n-1}=K_n/(K_{n,n-1})$ . Therefore

\begin{center}
    $[K_{n-1},K_{n-1} ] \cong [K_n,K_n]/([K_n,K_n]\cap K_{n,n-1})\cong (K_{n,n-1}\cdot[K_n,K_n])/K_{n,n-1}$
\end{center}
by the isomorphism theorems. \\

\noindent We substitute into the inductive step:
\begin{center}
    $(K_{n,n-1}\cdot[K_n,K_n])/K_{n,n-1} \subseteq (K_n/(K_{n,n-1}))^p
    \Rightarrow K_{n,n-1}\cdot[K_n,K_n]\subseteq (K_n/(K_{n,n-1}))^p \cdot K_{n,n-1}$
\end{center}

\noindent Since on the right hand side $K_{n,n-1}^p \cong \{ e\}$, and on the left hand side, $K_{n,n-1}$ is in the kernel of the surjection, we finally derive

\begin{center}
    $[K_n,K_n]\subset K_n^p \cdot K_{n,n-1}$.
\end{center}
 It is left to show that $K_{n,n-1}\subseteq K_n^p$. We show that $K_{n,n-1}\subset K_n^p$ by showing that for all $ h\in K_{n,n-1}$ there exists a $ g \in K_n$ such that $h=g^p.$ We need to show that for any $X\in M_{2 \times 2}(\Z/p^n) $ there exists a $ \ Y\in M_{2 \times 2}(\Z/p^n)$ s.t.
 \begin{center}
     $\1+p^{n-1}X \equiv \big( \1+pY)^p$ mod $p^n$,
 \end{center}
 where $\1+p^{n-1}X =h$ and $1+pY  =g$.\\
 
 \noindent Take $X= \begin{pmatrix} a & b \\ c & d \end{pmatrix}$, with $a,b,c,d \in \Z$.  Define further
\begin{center}$Y=p^{n-3}A,\ $ with $A= \begin{pmatrix} a & b \\ c & -a \end{pmatrix}$.
\end{center}

\noindent We need to verify that $\1+pY \in K_n$, i.e. that $\textnormal{det}(\1+pY)\equiv 1 \mod p$. Note that we assume $n \geq 3$ as $n=2$ was the base case in our inductive proof. Indeed, 
\begin{center}$ \textnormal{det}(\1+pY)=1+2p^{n-2}a-p^{2n-4}a^2-p^{2n-4}bc \equiv 1 \mod p.$ 
\end{center}Proceeding,
\begin{align*}
  (\1+pY)^p &= (\1+p(p^{n-3}A))^p \\
           &= (\1+p^{n-2}A)^p \\
           &= \1+p(p^{n-2}A)+\frac{p(p-1)}{2}(p^{n-2}A)^2+\{...\} \\
           &= \1+p^{n-1}A+p^n p^{n-3}A^2+\{...\}\\
           &\equiv \1+p^{n-1}A \mod p^n
\end{align*}
 
Substituting $A$:
 
 \begin{center} $(\1+p^{n-1}A) = \begin{pmatrix} 1+p^{n-1}a &p^{n-1}b \\ p^{n-1}c &1-p^{n-1}a \end{pmatrix}$
 \end{center}
 To show that $\1+p^{n-1}X \equiv (\1+pY)^p \mod p^n$ we show that $a \equiv -d \mod p$. This is a direct consequence of the fact that $det(X)=1$ implies $tr(X)=2$. Thus $-p^{n-1}a \equiv p^{n-1}d \mod p^n$ and hence $1+p^{n-1}X \equiv 1+p^{n-1}A \mod p^n$ as required.\\
 
\noindent We have thus finished the inductive proof that $[K_n,K_n]\subseteq K_n^p$. Hence $K_n$ is a powerful pro-$p$ group.\\

The proof for $L_n$ is analogous with relaxed conditions; we require the matrix to be only invertible. As we have shown the proof for determinant equal to 1 we have shown a stronger case. The case for the invertible matrix follows. \hfill $\square$\\

\begin{defn}
    For a group $G$ denote by $\Omega_1(G)$ the subgroup of $G$ generated by elements of order $p$. Then $G$ is called $\Omega$-extendable if $\Omega_1(G)$ is elementary abelian and there exists a central extension $E \rightarrow \tilde{G} \xrightarrow{\alpha} G$, where $E\cong \Omega_1(G)$ and every non-trivial element in $\Omega_1(G)$ is the image of an element of $\tilde{G}$ of order $p^2$.\\
\end{defn}

\begin{lem}\label{extendable} For odd primes $K_n$ and $L_n$ are $\Omega$-extendable for all $n>1$. 
\end{lem}
\textbf{Proof:} We consider $K_n$. Let $A\in K_{n,n-1} \triangleleft \ K_n $ be non-trivial. Then by Proposition \ref{elem.ab.} $A$ has order $p$. We need to construct $\tilde{K_n}$ such that for any such $A$ there exists $B\in \tilde{K_n}$ with $\alpha(B)=A$ and order of $B$ equal to $p^2$.\\

\noindent Let $E= K_{n+1,n}$. By Proposition \ref{elem.ab.} $E\cong (\Z/p)^3$. Then $E\cong K_{n,n-1}= \Omega_1(K_n)$. Let $\tilde{G}=K_{n+1}$.\\
This ensures that
\begin{center}
    $1 \rightarrow K_{n+1,n} \rightarrow K_{n+1} \rightarrow K_n \rightarrow 1$
\end{center}
is a short exact sequence.

\noindent We have the following diagram:\\

\begin{center}
\begin{tikzcd}
K_{n+1,n} \arrow[r] \ar[equal]{d}                                       & K_{n+1} \arrow[r]                                          & K_n                                              \\
K_{n+1,n} \arrow[r]         & K_{n+1,n-1} \arrow[r, "\gamma_{n+1,n}"] \arrow[u, hook]                      & K_{n,n-1} \arrow[u, hook]                        \\

\end{tikzcd}
\end{center}

Let \[
X= \begin{pmatrix} a & b \\ c & d \end{pmatrix},
\] $a,b,c,d \in \Z/p^{n}$, not all zero be such that $A=\1+p^{n-1}X \in K_{n,n-1}$.\\
We need $B= \1+p^{n-1}Y \in K_{n+1,n-1}$ such that $A=\gamma_{n+1,n} (B)$, $B^p$ is non-trivial, and $B^{p^2} \equiv \1 \mod p^{n+1}$. \\
We lift $a,b,c,d$ to $a',b',c',d'$ in $\Z/p^{n+1}$ where $a=a',\ b=a',\ c=c',\ d=d'$, creating the lift $X'=X$. Let $X'=X=Y$, then $B=\1+p^{n-1}X$. It is clear that $A=\gamma_{n+1,n} (B)$. We further find that $B^p$ is not trivial in $K_{n+1,n-1}$:
\begin{align*}
    B^p=(\1+p^{n-1}X)^p&=\sum_{k=0}^p \binom{p}{k}(p^{n-1}X)^k 
\end{align*}
Expanding the binomial coefficient gives
\begin{align*}
    & \frac{(2)(3)\cdot \cdot \cdot (p-k)(p-k+1) \cdot \cdot \cdot (p-1)(p)}{k!(2)(3) \cdot \cdot \cdot (p-k-1)(p-k)}p^{(n-1)k}X^k \\
    &=\frac{(p-k+1) \cdot \cdot \cdot (p-1)(p)}{k!}p^{nk-k}X^k
\end{align*}

For  $k=1$ this expression equals $\frac{p(p^{n-1})}{1}X=p^nX$. It follows that
\begin{center}
    $B^p=(\1+p^{n-1}X)^p \not\equiv \1 \mod p^{n+1}$.
\end{center}
Hence $B$ is non-trivial. Finally, we check that $|B|=p^2$, i.e. that $B^{p^2}\equiv \1 \mod p^{n+1}$.
\begin{align*}
    (\1+p^{n-1}X)^{p^2}&=\sum_{k=0}^{p^2} \binom{p^2}{k}(p^{n-1}X)^k 
\end{align*}
Expanding the binomial coefficient gives
\begin{align*}
    & \frac{(2)(3)\cdot \cdot \cdot (p^2-k)(p^2-k+1) \cdot \cdot \cdot (p^2-1)(p^2)}{k!(2)(3) \cdot \cdot \cdot (p^2-k-1)(p^2-k)}p^{(n-1)k}X^k \\
    &=\frac{(p^2-k+1) \cdot \cdot \cdot (p^2-1)(p^2)}{k!}p^{nk-k}X^k
\end{align*}
The expression $\binom{p^2}{k}p^{nk-k}$ is divisible by $p^{3+nk-2k}$. The exponent $3+nk-2k$ is greater than $n+1$ and thus $\binom{p^2}{k}p^{nk-k}$ is congruent to zero mod $p^{n+1}$ and thus
\begin{align*}
    (\1+p^{n-1}X)^{p^2} &=\1+\sum_{k=1}^{p^2} \binom{p^2}{k}(p^{n-1}X)^k \\
    &\equiv \1 \mod p^{n+1}
\end{align*}

We have confirmed the conditions laid out in Definition \ref{extendable}. This concludes that $K_n$ is $\Omega$-extendable. The proof for $L_n$ is analogous. \hfill $\square$ \\

\vspace{0.5cm}

In order to state and prove the main result of this section we need the following corollary and notation from [MS]:\\

\noindent \textbf{Notation:} For a $p$-group $G$ denote by $d(G)$ the minimal number of generators of $G$, and by $\Phi(G)=G^p[G,G]$ the Frattini subgroup of $G$.

\begin{thm}\label{mis} Let G be a powerful pro-p-group with $d(G)=m$ and $d(\Phi(G))=d$. 
Set $\Omega=\Omega_1(G)$ and $k=d(\Omega)$. The following are equivalent:\\

$(i)$  $\Omega$ is abelian and there exist $y_1,...,y_{k+d-m}$ in $H^2(G)$ and a basis $x_1,...,x_m$ of $H^1(G)$ such that\\

    $H^*(G)= \begin{cases} \Lambda[x_1,...,x_m] \otimes \Z/p[y_1....,y_{k+d-m},\beta x_{d+1},...,\beta x_m] \ \ $ for $ p>2\\ \Lambda[x_1,...,x_m] \otimes \Z/2[y_1....,y_{k+d-m}, x_{d+1},..., x_m] \ \ \ \ \ \ $ for $ p=2

\end{cases}$.\\

$(ii)$ $G$ is $\Omega$-extendable.
\end{thm}

\textbf{Proof:} See Page 239 in [MS]. \hfill $\square$\\

We are now in the position to state the following theorem:

\begin{thm}\label{Cohom(Kn)}
Let $K_n$ and $L_n$ be defined as before with $p>2$. 

\begin{center}
  $H^*(L_n,\mathbb{F}_p)=  \Lambda[x_1,x_2,x_3,x_4] \otimes \Z/p[y_1,y_2,y_3,y_4],  $\\
   $H^*(K_n,\mathbb{F}_p)=  \Lambda[x_1,x_2,x_3] \otimes \Z/p[y_1,y_2,y_3] $
\end{center}

 for all $n>1$, where $|x_i|=1,\ |y_i|=2$.
\end{thm}

\textbf{Proof:} 
By Proposition \ref{powerful} $K_n$ and $L_n$ are powerful pro-$p$-groups.\\
We first focus on $K_n$: $\Omega = \Omega_1(K_n)$ is  generated by the elements of $K_n$ order $p$. We use $\Omega \cong K_{n/n-1}$. In Proposition \ref{elem.ab.} we have shown that every element in $K_{n/n-1}$ has order $p$. This gives $m=d=k=3$. \\
For $L_n$ we have $\Omega = \Omega_1(L_n)$, $\Omega \cong K_{n/n-1}$, and $m=d=k=4$.\\
As shown in Lemma \ref{extendable} both $K_n$ and $L_n$ are $\Omega$-extendable. The statement of the theorem follows.\hfill $\square$\\

\section{The LHS spectral sequence}

The Lyndon-Hochschild-Serre spectral sequence is a well-established tool. We recall its definition:
\begin{defn}\label{LHS spsq}[Lyndon-Hochschild-Serre spectral sequence] 
Let $1 \rightarrow A \rightarrow B \rightarrow C \rightarrow 1$ be a group extension and $k$ a field. The Lyndon-Hochschild-Serre (LHS) spectral sequence associated to the extension has second page 
\begin{center}$E_2^{i,j}=H^i(C,H^j(A,k)) \Rightarrow H^{i+j}(B,k),$
\end{center}
meaning that the bigraded modules $E_2^{i,j}$ converge to $H^{i,j}$ via the differentials $d_r$ and the lifting of the $E_{\infty}$ page.
\end{defn}

\noindent We apply the LHS spectral sequence to the group extensions

\begin{center}
    $ 1 \rightarrow L_n \rightarrow S_p(n,GL) \rightarrow S_p(1,GL)\cong C_p \rightarrow 1 \hspace{2cm}$ $(*)$\\
    $ 1 \rightarrow K_n \rightarrow S_p(n,SL) \rightarrow S_p(1,SL) \cong C_p \rightarrow 1\hspace{2cm}$ $(**)$
\end{center}

\noindent \textbf{Notation:} We denote by $E_r(S_p(n,GL))$ and $E_r(S_p(n,SL))$ the $r^{th}$ pages of the LHS spectral sequences associated to the above group extensions.\\

\noindent In order to compute the cohomology groups of the Sylow $p$-subgroups we need the cohomologies of $L_n$ and $K_n$ respectively, as well as the cohomology of $C_p$. It is well-known that for odd primes $p$
\begin{center}
    $H^*(C_p, \mathbb{F}_p) \cong \Lambda(\zeta_1) \otimes \mathbb{F}_p(\zeta_2)$ with $|\zeta_1|=1, \ |\zeta_2|=2$ 
\end{center}
\vspace{0.5cm}

We examine the cohomology rings of $L_n$ and $K_n$. The following is one of the main theorems of this paper.

\begin{thm}\label{G-module isomorphism} For $p$ an odd prime denote $M:=H^1(K_n, \F_p)$ and $N:=H^1(L_n, F_p)$. By $(*)$ and $(**)$ these have a $C_p$ action. In fact, for any $n>1$
\begin{center}
$H^*(L_n,\F_p)=S(N) \otimes \Lambda(N)$\\
    $H^*(K_n,\F_p)=S(M) \otimes \Lambda(M)$
    
\end{center}
as $C_p$-modules, where $\Lambda$ denotes the exterior algebra and $S$ the symmetric algebra.
\end{thm}

\textbf{Proof:}
We focus on $H^*(K_n,\F_p)$. By Theorem \ref{Cohom(Kn)} 
\begin{center}
    $H^*(K_2,\mathbb{F}_p) \cong H^*(K_n,\mathbb{F}_p) \cong \Lambda[x_1,x_2,x_3] \otimes S[y_1,y_2,y_3] $
\end{center} as rings with $|x_i|=1,\ |y_i|=2$.  
By Definition \ref{LHS spsq}, the $E_2$-page of the LHS spectral sequence associated to the extension $(**)$ is given by
\begin{center}$E_2^{i,j}=H^i(C_p,H^j(K_n,\F_p)) \Rightarrow H^{i+j}(Syl_p(n,S),\F_p)$,
\end{center}
making $H^*(K_n,\F_p)$ a $C_p$ module.
\noindent We start with $H^1(K_n,\F_p)$.
By definition, $H^1(K_n,\F_p) \cong Hom_{\F_p}(K_n, \F_p)$. 
For $\alpha \in Hom_{\F_p}(K_n, \F_p)$, $g \in C_p$  and $x \in K_n$,
\begin{center}
    $g\cdot \alpha(x)=\alpha(gxg^{-1})$,
\end{center}
 thus $C_p$ acts on $K_n$ by conjugation. We denote this action by $\rho_{K_n}$. We show that for any $n>2$, $\rho_{K_n}$ and $\rho_{K_2}$ give the same action of $C_p$ on the cohomology. By the short exact sequence
 \begin{center}
    $1 \rightarrow K_{n,n-1}\rightarrow K_n \rightarrow K_{n-1}\rightarrow 1$,
\end{center}
$K_{n,n-1}$ is a normal subgroup of $K_n$, and hence fixed by conjugation. Thus $\rho_{K_n}$ restricts to $\rho_{K_{n,n-1}}$. By Proposition \ref{elem.ab.} $K_{2,1} \cong K_{n+1,n}$ via the map $\theta$ 
\begin{center}
    $\theta: K_{2,1}\rightarrow K_{n+1,n}$\\
    $\1+pX \mapsto \1+ p^nX$,
\end{center}
 and hence on the level of cohomology the isomorphism is compatible with the $C_p$-action. Thus  $H^1(K_{2},\mathbb{F}_p) \cong H^1(K_n,\mathbb{F}_p)$ as $C_p$-modules.\\

\noindent We move to $H^2(K_n,\F_p)$. This group as $\F_p$-vector space has dimension 6. We wish to determine the subspaces and examine the $C_p$-action on them. 

\noindent As $K_n$ is a powerful group, $\Lambda^2(H^1(K_n,\F_p))$ injects into $H^2(K_n,\F_p)$, which is the statement of Theorem 5.1.6 in [SW] and has been proven there, creating a subspace of $H^2(K_n,\F_p)$. This exterior algebra has dimension 3 and inherits the $C_p$-action from $H^1(K_n)\cong \Lambda^1(x_1,x_2,x_3)$.\\

\noindent In order to examine the other subspace of $H^2(K_n,\F_p)$, we turn our attention to the short exact sequence
\begin{center}
    $1 \rightarrow K_{n+1,n} \rightarrow K_{n+1} \rightarrow K_n \rightarrow 1$
\end{center}
 There is a LHS spectral sequence associated to it, with second page given by
 \begin{center}
    $E_2^{i,j}=H^i(K_n,H^j(K_{n+1,n}\F_p))$.
\end{center}
The $E_2$ page has bottom left corner

\begin{center}
\begin{tikzcd}
{}                          &                                     &          &          &    \\
                            & ({H^1(K_{n+1,n})})^{K_n} \arrow[rrd, "\delta"] &          &          &    \\
                            & 1                                   & H^1(K_n) & H^2(K_n) &    \\
{} \arrow[uuu] \arrow[rrrr] &                                     &          &          & {}
\end{tikzcd}
\end{center}

\noindent 
 We have shown above that the isomorphism $K_{2,1} \cong K_{n+1,n}$ is compatible with the $C_p$-action and hence $H^1(K_{2,1},\mathbb{F}_p) \cong H^1(K_{n+1,n},\mathbb{F}_p)$ as $C_p$-module. The map $\delta$ is injective, which is seen as follows: Assume $\delta$ is not an injection, then $\textnormal{ker}(\delta)\neq 0$. This means that $E_3^{0,1}\neq 0$, and since all higher differentials vanish due to slope and length of the arrows, $E_{\infty}^{0,1}\neq 0$. Consequently $H^1(K_{n+1},\mathbb{F}_p)$ consists of more than $H^1(K_n,\mathbb{F}_p)$. But $H^1(K_{n},\mathbb{F}_p) \cong H^1(K_{n+1},\mathbb{F}_p)$ as groups as stated earlier in this proof. Hence $E_{\infty}^{0,1}=0$ and in particular $E_3^{0,1}= 0$ due to length and slopes of the higher differentials. This can only happen when $\textnormal{ker}(\delta)=0$ and hence $\delta$ is an injection.\\

\noindent We have now the following information about $H^2(K_n,\F_p)$: As vector space it has dimension 6, one subspace is $\Lambda^2(H^1(K_n,\F_p))$, one subspace is $S^1(y_1,y_2,y_3)$, and one subspace is $\textnormal{im}(\delta)$. Each of these subspaces has dimension 3. Looking at the $E_3$-page we have that $E_3^{2,0}\cong H^2(K_n)/\textnormal{im}(\delta)$. As $K_{n+1}$ is powerful, $\Lambda^2(H^1(K_n,\F_p))$ must survive onto the $E_{\infty}$ page, and thus $\Lambda^2(H^1(K_n,\F_p)) \cap \textnormal{im}(\delta)=0.$ This gives that $H^2(K_n,\F_p)\cong \Lambda^2(H^1(K_n,\F_p)) \oplus \textnormal{im}(\delta)$ as $\F_p$-vector spaces, and thus $S^1(y_1,y_2,y_3)\cong \textnormal{im}(\delta)$. As $C_p\cong (\F_p)$, and $\delta$ is an injective map of $C_p$-modules, $S^1(y_1,y_2,y_3)$ is a $C_p$-module.\\

\noindent We have thus proven the statement of the theorem for $K_n$. The proof for $L_n$ is analogous. \hfill $\square$ \\

This theorem gives the following corollary:

\begin{cor}\label{E_2 all p}
With the previously established notation there are group isomorphisms for $p$ an odd prime
\begin{center}
 $E_2^{i,j}(S_p(2,GL))\cong E_2^{i,j}(S_p(n,GL)) \ \textnormal{ for all } n>2$\\
 $E_2^{i,j}(S_p(2,SL))\cong E_2^{i,j}(S_p(n,SL)) \  \textnormal{ for all }  n>2$.
 \end{center}
\end{cor}

\textbf{Proof:}
Follows immediately from Theorem \ref{G-module isomorphism} and Definition \ref{LHS spsq}. \hfill $\square$\\

Note that this corollary gives information only about the $E_2$-pages of the LHS-spectral sequences associated to the presented extens of the Sylow $p$-subgroups. It does not make claims about the convergence of the spectral sequence, or on the lifting of generators and relations. These questions need to be answered computationally.

\vspace{0.5cm}

 \section{Stable Elements}

 The motivation to study the cohomology of the Sylow $p$-subgroups of $GL_2(\Z/p^n)$ and $SL_2(\Z/p^n)$ was that the stable elements of these cohomologies give the modular cohomologies of $GL_2(\Z/p^n)$ and $SL_2(\Z/p^n)$. In this section we give an explicit description of what the stable elements look like. For this we will use the toolbox provided by the theory of fusion systems. \\
 Recall that for $K\subset G$ and $g\in G$, $c_g:H^*(K)\rightarrow H^*(^gK)$ is the map induced by conjugation by $g$. Further recall the restriction map on cohomology, $res^G_K :H^*(G)\rightarrow H^*(K)$. In their 1956 book \textit{Homological Algebra} ([CE]) Cartan and Eilenberg give the following definition of stable elements:
\begin{defn}
  Let $G$ be a finite group, $S$ a subgroup, $A$ a $G$-module. For $g\in G$,  $a\in H^*(S,A)$ and \\
  $ c_ga\in H^*(^gS,A)$, $a$ is called stable if
\begin{center}
$res^S_{S\cap ^gS}(a)=res^{^gS}_{S\cap ^gS}(c_ga)$.
\end{center}
\end{defn}
 When moving to $S$ a Sylow $p$-subgroup of $G$ we have the following desired result:
\begin{lem}\label{stable elem method}
    Let $S$ be a Sylow $p$-subgroup of a finite group $G$. Denote by $\Theta$ the set of stable elements in $H^*(S,\F_p)$. There is a ring isopmorphism
        $\Theta \cong H^*(G,\F_p)$.
   
\end{lem}

\noindent \textbf{Proof:} \\
We have the following diagram 
\begin{center}
\begin{tikzcd}
                                                                  
H^i(G) \arrow[rr, "res_S^G"] \arrow[rrrr, "\times|G:S|"', bend right] &  & H^i(S) \arrow[rr, "tr^G_S", hook] &  & H^i(G)
\end{tikzcd}
\end{center}
where $res_S^G$ and $tr^G_S$ are the usual restriction and transfer maps. By Theorem 10.1 in [CE],  $\textnormal{im}(res_S^G) = \Theta$. Furthermore, by Proposition 6.5 in [AM], if $a\in H^*(S)$ is stable, then $  res_S^G\circ tr^G_S(a)=|G:S|a$.
 Since $S$ is a Sylow $p$-subgroup of G, $|G:S|$ is co-prime to $p$, and thus the composition $ res_S^G\circ tr^G_S $ is injective. This combined with $\textnormal{im}(res_S^G) = \Theta$ shows that $\Theta$ is mapped by $tr_S^G$ injectively into $H^*(G)$. By Theorem 10.1 in [CE] this mapping is also surjective. \\
 Thus $H^*(G)\cong \Theta$ as claimed. \hfill $\square$\\

Unfortunately, as pointed out in section \ref{section Kn}, classical group theory does not get us far. If we let $G=GL_2(\Z/p^n)$ for example, the order of $G$ grows very large very quickly as we increase $n$, supplying a rapidly growing number of elements to conjugate with. Furthermore, our Sylow $p$-subgroups are neither abelian nor normal in $G$, making a direct application of the theory as set out in [CE] difficult. In order to describe the stable elements in $H^*(S_p(n),\F_p)$ regardless, we turn to the theory of fusion systems. As a deeper background read we suggest to consult [AKO] and [BLO]. We continue by recalling the required notions from the theory of fusion systems. The following is Definition 1.1 in [BLO]. \\

\begin{defn}[Fusion system] A fusion system $\fs$ on a finite p-group S is a category with objects given by the subgroups of $S$ and morphisms form the set $Hom_{\fs}(P,Q)$ satisfying \\
i) $Hom_S(P,Q)\subset Hom_{\fs}(P,Q)\subset Inj(P,Q)$ for all $P,Q \leq S$\\
ii) Every morphism in $\fs$ factors through an isomorphism in $\fs$ followed by an inclusion.
\end{defn}

We can further impose two technical axioms resembling the Sylow axioms in classical group theory (consult [BLO] Definition 1.2 for further information). Fusion systems satisfying these are called \textit{saturated fusion systems}.\\
 Let Take $G$ a finite group, $S\in Syl_p(G)$. It is natural to think of a fusion system attached to this, $\fs=\fs_S(G)$. We call such a fusion system a \textit{group fusion system}, and all morphisms are induced by conjugation with elements $g \in G$. The study of fusion systems often focusses only on a $p$-group $S$. In the case of so-called \textit{exotic fusion systems}, the fusion systems that cannot even be realised by a group $G$, meaning $\fs$ has morphism that do not come from conjugation with elements in a larger group. Both kinds of fusion systems are saturated. For further information we recommend [AKO] and [BLO] and [BLO2]. In our case, however, we do not study fusion systems in their full generality, but rather use the toolbox the field provides to understand the cohomology of the groups we are interested in. The following theorem provides the basis for this approach, as it generalises the stable elements theorem from Cartan and Eilenberg.

\begin{thm}
    For a finite group $G$ with Sylow $p$-subgroup $S$,
    \begin{center}
        $H^*(G,\F_p)\cong H^*(\mathcal{F}_S(G),\F_p)$.
        \end{center}
\end{thm}
\textbf{Proof:} Follows from Proposition 1.1 in [BLO2].\\

\noindent In order to state the stable elements in the fusion systems language we need the following definitions:

\begin{defn}[$\mathcal{F}$-centric $p$-group]
    Let $\mathcal{F}$ be a fusion system over a $p$-group $S$. A subgroup $P\leq S$ is $\mathcal{F}$-centric if $P$ and all of its $\mathcal{F}$-conjugates contain their $S$-centralisers.
\end{defn}

\begin{defn}
    Denote by $\mathcal{O}(\fs)$ the orbit category of $\fs$ over a p-group S. This category has objects which are the subgroups of S whose morphisms are defined by
    \begin{center}
        $Mor_{\mathcal{O}(\fs}(P,Q)=Rep_{\fs}(P,Q)=Inn(Q)\backslash Hom_{\fs}(P,Q)$.
    \end{center}
    We let $\mathcal{O}^c(\fs)$ denote the full subcategory of $\mathcal{O}(\fs)$ whose objects are the $\fs$-centric subgroups of $\fs$.
\end{defn}

The stable elements formulation in the setting of fusion systems is as follows:
 \begin{defn}
 For a Sylow $p$-subgroup $S$ of $G$ and $\alpha\in H^*(S)$, $(res^S_P(\alpha))_{P\in \mathcal{F}^c(S)} \in \underset{\xleftarrow[\mathcal{O}^c(\mathcal{F})]{}}{\mathrm{lim}} H^*(-)$ if and only if
 \begin{center}
     $\forall \phi :P\rightarrow P' \ \in Hom_{\mathcal{O}^c(\mathcal{F})}(P,P')$ \\
     $\phi^*(res^S_{P'}(\alpha))=res^S_P(\alpha)$.
 \end{center}
 \end{defn}

Note that this inverse limit runs over all $\fs$-centric subgroups of $S$. This can be a large number, and cumbersome to compute. Hence we would like to restrict our attention. This is possible as shown after the next definition.

\begin{defn}[p-radical subgroup]
For a finite group $G$ with subgroup $P$, $P$ is called p-radical if $N_G(P)/P$ has no non-trivial normal p-subgroups, i.e. $O_p(N_G(P)/P)=1$, where $O_p(-)$ denotes the maximal normal p-subgroup. 
\end{defn}

\begin{prop}\label{centric radical cohom} Let $\fs=\fs_G(S)$ be a group fusion system over a Sylow $p$-subgroup $S$. In order to compute the modular cohomology of a fusion system it is sufficient to take the $\fs$-centric and $\fs$-radical subgroups of $S$ into consideration.
\end{prop}

Before we can give the proof we need the following definitions, which are Definition 1.7 and 1.8 in [BLO].
\begin{defn}
    Let $\fs$ be a fusion system over a p-group S. A centric linking system associated to $\fs$ is a category $\mathcal{L}$ whose objects are the $\fs$-centric subgroups of S togethr with a functor
    \begin{center}
        $\pi:\mathcal{L} \rightarrow \fs^c$,
    \end{center}
    and "distinguished" monomorphisms $P \xrightarrow{\delta_P} Aut_{\mathcal{L}}(P)$ for each $\fs$-centric subgroup $P\leq S$, which satisfy the following conditions.\\
(A) $\pi$ is the identity on objects and surjective on morphisms.  More precisely, for each pair of objects $P,Q \in \mathcal{L}$, $Z(P)$ acts freely on $Mor_{\mathcal{L}}$ by composition (upon identifying $Z(P)$ with $\delta_P(Z(P)\leq Aut_{\mathcal{L}}(P))$, and $\pi$ induces a bijection
\begin{center}
    $Mor_{\mathcal{L}}(P,Q)/Z(P) \xrightarrow{\cong} Hom_{\fs}(P,Q).$
\end{center}
(B) For each $\fs$-centric subgroup $P\leq S$ and each $g\in P$, $\pi$ sends $\delta_P(g)\in Aut_{\mathcal{L}}(P)$ to $c_g\in Aut_{\fs}(P)$.\\
(C)For each $f\in Mor_{\mathcal{L}}(P,Q)$ and each $g \in P$ the following square commutes in $\mathcal{L}$:
\begin{center}
\begin{tikzcd}
P \arrow[r, "f"] \arrow[d, "\delta_P(g)"] & Q \arrow[d, "\delta_Q(\pi(f)(g))"] \\
P \arrow[r, "f"]                          & Q                                 
\end{tikzcd}.
\end{center}
\end{defn}

\begin{defn}
    A p-local finite group $(S,\fs, \mathcal{L})$, where $\fs$ is a saturated fusion system over the p-group S and $\mathcal{L}$ is a centric linking system associated to $\fs$. The classifying space of the p-local group $(S,\fs,\mathcal{L})$ is the space $|\mathcal{L}|^{\land}_p$.
\end{defn}
 \noindent We denote by $\mathcal{L}_0\subseteq \mathcal{L}$ the full subcategory which contains all $\fs$-radical and $\fs$-centric subgroups of $S$.\\

\noindent We can now provide the proof for Proposition \ref{centric radical cohom}.\\

\textbf{Proof:} Following the results of [BLO], we prove the following sequence of homotopy equivalences and cohomology isomorphisms:
\begin{center}
    $H^*(|L_0|^{\wedge}_p,\F_p) \ \cong_{(i)} \ H^*(|L|^{\wedge}_p,\F_p) \ \cong_{(ii)} \ H^*(|L|,\F_p) \ \cong_{(iii)} \ H^*(\fs,\F_p) \ \cong_{(iv)} \ H^*(G,\F_p)$.
\end{center}
The first isomorphism is a consequence of Corollary 3.6 in [BLO], which states that there is a mod p homotopy equivalence $|L_0|^{\wedge}_p \simeq |L|^{\wedge}_p$. Homotopy equivalence induces isomorphism on cohomology, and our cohomology coefficients are $\F_p$. Isomorphism $(ii)$ follows from the fact that we have $\F_p$ as cohomology coefficients. For isomorphism $(iii)$ we apply Theorem 5.8 in [BLO]. Isomorphism $(iv)$ has been stated above.

\qed\\

For our case where $S$ is a Sylow $p$-subgroup of $SL_2(\Z/p^n)$ or $GL_2(\Z/p^n)$ we still run into the problem that the number of subgroups of $S$ grows very quickly as $n$ increases. For example, a Sylow 3-subgroup $SL_2(\Z/3^2)$ has 20 conjugacy classes of subgroups. For a Sylow 3-subgroup of $SL_2(\Z/3^3)$ this number is 97, and for a Sylow 3-subgroup of $SL_2(\Z/3^4)$ we are at 282. To require examination of every single subgroup to determine whether it is $\fs$-centric and $p$-radical would make our computations impossible. The proposition following the next definition prevents this.

\begin{defn}
    Take $\fs$ a fusion system over a p-group S, and fix $Q \leq S$. Q is normal in $\fs$, write $Q \triangleleft \fs$, if $Q \triangleleft S$, and for all $P,R \leq S$ and
for all $\phi \in Hom_{\fs}(P,R)$, $\phi$ extends to a morphism $\tilde{\phi}\in Hom_{\fs}(PQ,RQ)$ such that $\tilde{\phi}(Q)=Q$.
\end{defn}

\noindent Note that in our case $K_n$ is not only normal in $S_p(n,SL)$ but also normal in $SL_2(\Z/p^n)$. Since the fusion system on the Sylow $p$-subgroup $\fs_{S_p(n,SL)}(SL_2(\Z/p^n))$ is a group fusion system, all maps in the fusion system come from conjugation by elements in $G$. As $K_n$ is normal in $SL_2(\Z/p^n)$, $K_n$ is also normal in $\fs_{S_p(n,SL)}(SL_2(\Z/p^n))$. The same applies for $L_n \triangleleft \fs_{S_p(n,GL)}(GL_2(\Z/p^n))$.

\begin{prop}\label{Craven}
    If $\mathcal{F}$ is a saturated fusion system over a p-group S and $Q \triangleleft \mathcal{F}$ then $Q\leq R$ for every $\mathcal{F}$-centric and $p$-radical $R\leq S$.
\end{prop}
\textbf{Proof:} Page 118 in [C] \qed\\

We are now in the position to state and prove the main theorem of this paper. Recall that by Lemma \ref{stable elem method} the stable elements of the modular cohomology a Sylow $p$-subgroup of $G$ give us the modular cohomology of $G$. The following theorem thus describes the modular cohomology rings of $SL_2(\Z/p^n)$ and $GL_2(\Z/p^n)$ for odd $p$.

\begin{thm}\label{final theorem} Let $n>1$ and $p$ be an odd prime.
    Denote by $\Theta_{n}$ ring of stable elements  in $H^*(S,\F_p)$, where $S$ is a Sylow $p$-subgroup of $SL_2(\Z/p^n)$. Then
    \begin{center}
        $\Theta_n= \{H^*(S,\F_p)^{N_{SL_2(\Z/p^n)}(S)/S}\} \cap \{x\in H^*(S,\F_p) \ | \  res^{S}_{K_n}(x)\in H^*(K_n,\F_p)^{N_{SL_2(\Z/p^n)}(K_n)/K_n} \}$.
    \end{center}
    Denote by $\Psi_{n}$ ring of stable elements  in $H^*(Q,\F_p)$,  where $Q$ is a Sylow $p$-subgroup of $GL_2(\Z/p^n)$. Then
    \begin{center}
        $\Psi_n= \{H^*(Q,\F_p)^{N_{SL_2(\Z/p^n)}(Q)/Q}\} \cap \{x\in H^*(Q,\F_p) \ | \  res^{Q}_{L_n}(x)\in H^*(L_n,\F_p)^{N_{SL_2(\Z/p^n)}(L_n)/L_n} \}$.
    \end{center}
    
\end{thm}

\textbf{Proof:} We show the proof for $SL_2(\Z/p^n)$. As seen above for a fusion system $\mathcal{F}$ over a $p$-group S, $H^*(\mathcal{F})=\underset{\xleftarrow[\mathcal{O}^c(\mathcal{F})]{}}{\mathrm{lim}} H^*(-)$, where the limit runs over all $\mathcal{F}$-centric, $p$-radical subgroups $P$. We are thus tasked to determine such $P$ and their outer automorphism groups. For a $p$-subgroup $P$ of a group $G$ the outer automorphism group is $\textnormal{Out}_{G}(P)=\textnormal{Aut}_G(P)/\textnormal{Inn}_G(P)=N_G(P)/C_G(P)\cdot P$. Note that if $P$ is $\mathcal{F}$-centric then $C_G(P)\leq P$ and thus $\textnormal{Out}_{G}(P)=N_G(P)/ P$.\\
Clearly $S$ is $\mathcal{F}$-centric. To show that is is $p$-radical we need its normaliser. As $S$ is a Sylow $p$-subgroup, it has maximal $p$-rank, and hence the normaliser quotiented out by $S$ does not contain non-trivial normal $p$-groups. 
Another $\mathcal{F}$-centric subgroup is $K_n$. A computation shows that $S$ is not in the centraliser of $K_n$, and as $K_n$ is maximal and normal in $S$, this means that $K_n$ contains $C_{S}(K_n)$. 
In addition to $K_n$ being $\mathcal{F}$-centric we show that $K_n$ is p-radical: $N_G(P)/P = SL_2(\Z/p^n)/K_n=SL_2(\Z/p)$. This has $C_p$ as p-subgroup, but $C_p$ is not normal in $SL_2(\Z/p)$ and hence $K_n$ is p-radical.\\ 
 Thus, the $K_n$-component adds the second part of the inverse limit, with outer automorphism group $SL_2(\Z/p)$.\

\noindent We continue to show that no other subgroups of $S(n)$ are both $\mathcal{F}$-centric and $p$-radical.  By Proposition \ref{Craven} any normal subgroup of $\fs$ is included in every $\mathcal{F}$-centric and $p$-radical subgroup of $S$. As $K_n$ is normal in $\fs$, it follows that there cannot be any other $\mathcal{F}$-centric and $p$-radical subgroups of $S$ but the ones we already listed. \\
The proof for $GL_2(\Z/p^n)$ is analogous. \qed\\

\textbf{Remark:} To close this paper we would like to draw attention to limitations of this theorem. Whilst it gives a formula for computing the desired cohomology rings, we cannot claim that 
\begin{center}
    $H^*(GL_2(\Z/p^n),\F_p) \cong H^*(GL_2(\Z/p^{n+1}),\F_p)$.
\end{center}
The reason for this is that we do not know whether the cohomology of the Sylow $p$-subgroups is the same across all $n>1$. This has been remarked following Corollary \ref{E_2 all p}. It is a topic of ongoing research to compute these cohomologies explicitly for each prime.

\hfill

\section{Bibliography}

\noindent \textbf{[A]} J. AGUADE: \textit{ The cohomology of the GL 2 of a finite field}. Arch. Math 34, 509–516 (1980).

\noindent \textbf{[AM]} Alejandro ADEM, James MILGRAM: \textit{Cohomology of finite groups}. Springer Verlag Berlin Heidelberg, 2004 

\noindent \textbf{[BLO]} Carles BROTO, Ran LEVI, Bob OLIVER: \textit{The homotopy theory of fusion systems}. Journal of the American Mathematical Society 16 (4), 779-856  

\noindent \textbf{[BLO2]} Carles BROTO, Ran LEVI, Bob OLIVER: \textit{Homotopy equivalences of p-completed classifying spaces of finite groups}. Invent. math. 151, 611–664 (2003).

\noindent \textbf{[C]} David A. CRAVEN: \textit{Theory of Fusion Systems}. Cambridge University Press, 2012

\noindent \textbf{[CE]} Henri CARTAN, Samuel EILENBERG: \textit{Homological Algebra}. Princeton University Press, 1956

\noindent \textbf{[DGMP]} Antonio DIAZ, Oihana GARAIALDE, Nadia MAZZA, Sejong PARK: \textit{On the cohomolomology of pro-fusion systems}. Journal of Algebra and Its Applications 22 (11), 2350243

\noindent \textbf{[G]} K. GRUNBERG: \textit{Profinite Groups}. Algebraic Number Theory, Proceedings of an instructional conference, 116-127, Academic Press, 1967

\noindent \textbf{[K]} Bruno KAHN \textit{A characterization of powerfully embedded normal subgroups of a p-group}. J. Algebra 188 (1997), 401-408.

\noindent \textbf{[MS]} Pham Anh MINH, Peter SYMONDS: \textit{The Cohomology of Pro-p Groups with a Powerfully Embedded Subgroup}. Journal of Pure and Applied Algebra 189, 221-246, 2004.

\noindent \textbf{[Q]} Daniel QUILLEN: \textit{On the Cohomology and K-Theory of the General Linear Groups Over a Finite Field}
Ann. of Math. 96, 1972.

\noindent \textbf{[SW]} Peter SYMONDS, Thomas WEIGEL: \textit{Cohomology of $p$-adic Analytic Groups}. New Horizons in Pro-p Groups, du Sautoy, Segal and Shalev eds., Birkhauser (2000) 349-416 .\\

\end{document}